\documentclass[a4paper]{amsart}

\usepackage{amssymb,latexsym}
\usepackage{amsmath}
\usepackage{amsfonts}
\usepackage{amsthm}
\usepackage{amssymb}

\theoremstyle{plain}
\newtheorem{theorem}{Theorem}

\theoremstyle{definition}

\title[Rational Diophantine sextuples]
{Rational Diophantine sextuples containing two regular quadruples and one regular quintuple}
\begin{document}

\date{}


\author[A. Dujella]{Andrej Dujella}
\address{
Department of Mathematics\\
Faculty of Science\\
University of Zagreb\\
Bijeni{\v c}ka cesta 30, 10000 Zagreb, Croatia
}
\email[A. Dujella]{duje@math.hr}

\author[M. Kazalicki]{Matija Kazalicki}
\address{
Department of Mathematics\\
Faculty of Science\\
University of Zagreb\\
Bijeni{\v c}ka cesta 30, 10000 Zagreb, Croatia
}
\email[M. Kazalicki]{matija.kazalicki@math.hr}

\author[V. Petri\v{c}evi\'c]{Vinko Petri\v{c}evi\'c}
\address{
Department of Mathematics\\
Faculty of Science\\
University of Zagreb\\
Bijeni{\v c}ka cesta 30, 10000 Zagreb, Croatia
}
\email[V. Petri\v{c}evi\'c]{vpetrice@math.hr}

\begin{abstract}
A set of $m$ distinct nonzero rationals $\{a_1,a_2,\ldots,a_m\}$
such that  $a_ia_j+1$ is a perfect square for all $1\leq i<j\leq m$, is called a rational Diophantine
$m$-tuple. It is proved recently that there are infinitely many rational Diophantine sextuples.
In this paper, we construct infinite families of rational Diophantine sextuples with special structure,
namely the sextuples containing quadruples and quintuples of certain type.
\end{abstract}


\maketitle

\section{Introduction}

A Diophantine $m$-tuple is a set of $m$ distinct positive
integers with the property that the product of any two of its distinct
elements plus $1$ is a square. Fermat found the first Diophantine quadruple in integers $\{1,3,8,120\}$.
If a set of $m$ nonzero rationals
has the same property, then it is called
a rational Diophantine $m$-tuple.
The first example of a rational Diophantine quadruple was the set
$$
\left\{\frac{1}{16},\, \frac{33}{16},\, \frac{17}{4},\, \frac{105}{16}\right\}
$$
found by Diophantus. Euler proved that the exist
infinitely many rational Diophantine quintuples (see \cite{Hea}),
in particular he was able to extend the integer Diophantine quadruple found by Fermat,
to the rational quintuple
$$
\left\{ 1, 3, 8, 120, \frac{777480}{8288641} \right\}.
$$
Stoll \cite{Stoll} recently showed that this extension is unique.
Therefore, the Fermat set $\{1,3,8,120\}$ cannot be extended to a rational Diophantine sextuple.

In 1969, using linear forms in logarithms of algebraic numbers and a reduction method
based on continued fractions, Baker and Davenport \cite{B-D}
proved that if $d$ is a positive integer such that
$\{1, 3, 8, d\}$ forms a Diophantine quadruple, then $d$ has to be $120$.
This result motivated the conjecture that there does not exist a Diophantine quintuples in integers.
The conjecture has been proved recently by He, Togb\'e and Ziegler \cite{HTZ}
(see also \cite{BTF,duje-crelle}).

In the other hand, it is not known how large can be a rational Diophantine tuple.
In 1999, Gibbs found the first example of rational Diophantine sextuple \cite{Gibbs1}
$$
\left\{ \frac{11}{192}, \frac{35}{192}, \frac{155}{27}, \frac{512}{27}, \frac{1235}{48}, \frac{180873}{16} \right\}.
$$
In 2017 Dujella, Kazalicki, Miki\'c and Szikszai \cite{DKMS} proved that there are
infinitely many rational Diophantine sextuples, while
Dujella and Kazalicki \cite{Duje-Matija} (inspired by the work of Piezas \cite{P})
described another construction of parametric families of rational Diophantine sextuples.
Recently, Dujella, Kazalicki and Petri\v{c}evi\'c \cite{DKP-sext}
proved that there are infinitely many rational Diophantine sextuples
such that denominators of all the elements (in the lowest terms) in the sextuples are perfect squares.
No example of a rational Diophantine septuple is known.
The Lang conjecture on varieties of general type implies that
the number of elements of a rational Diophantine tuple is bounded by an absolute constant
(see Introduction of \cite{DKMS}).
For more information on Diophantine $m$-tuples see the survey article \cite{Duje-notices}.

Although the constructions of infinitely families of rational Diophantine sextuples
in \cite{Duje-Matija} and \cite{DKMS} are essentially different,
they have one common feature. Namely, in both constructions
(and also in \cite{DKP-sext}, which is a special case of \cite{Duje-Matija})
the sextuples contain two regular rational Diophantine quintuples.
The quintuple $\{a,b,c,d,e\}$ is called regular if
$$ (abcde+2abc+a+b+c-d-e)^2=4(ab+1)(ac+1)(bc+1)(de+1) $$
(see \cite{duje-quint,duje-irr,duje-sext,G2}). Similarly, the quadruple $\{a,b,c,d\}$ is called regular if
$$ (a+b-c-d)^2=4(ab+1)(cd+1) $$
(see \cite{duje-irr} for characterization of regular Diophantine quadruples and quintuples
in terms of corresponding elliptic curves).

In \cite{G3}, Gibbs collected over 1000 examples of rational Diophantine sextuples with relatively
small numerators and denominators. These examples are also sorted in \cite{G3} according to their
structure, which includes information of regular quadruples and quintuples which they contain.
We have extended the search for sextuples with small height and included also examples
with mixed signs (in \cite{G3} only sextuples with positive elements were considered).
We have observed a significant number of sextuples which contain
exactly one regular Diophantine quintuple and two regular Diophantine quadruples.
Thus, in this paper we study rational Diophantine sextuples having this structure.
Our main result is the following theorem.

\begin{theorem} \label{tm1}
There are infinitely many rational Diophantine sextuples which contain one regular Diophantine quintuple
and two regular Diophantine quadruples.
\end{theorem}

\section{Parametrizations of Diophantine triples} \label{sec:triples}

Let $\{a_1,a_2,a_3\}$ be a rational Diophantine triple and let $a_2=(r^2-1)/a_1$ and $a_3=(s^2-1)/a_1$
for rationals $r$ and $s$. By putting $a_2a_3+1=(a_2s^2-a_2+a_1)/a_1 = (1+(s-1)t)^2$, we get
$$ a_3=\frac{-4t(t-1)(a_1t-a_2)}{(-a_2+a_1t^2)^2}. $$
This parametrization of Diophantine triples was used in \cite{duje-sext} in construction of certain
rational Diophantine sextuples. Here we will use an equivalent, but simpler and more aesthetic parametrization
due to Lasi\'c \cite{luka}, which is symmetric in the three involved parameters:
\begin{align*}
a_1 &=  \frac{2 t_1 (1 + t_1 t_2 (1 + t_2 t_3))} {(-1 + t_1 t_2 t_3) (1 + t_1 t_2 t_3)}, \\
a_2 &=  \frac{2 t_2 (1 + t_2 t_3 (1 + t_3 t_1))} {(-1 + t_1 t_2 t_3) (1 + t_1 t_2 t_3)}, \\
a_3 &=  \frac{2 t_3 (1 + t_3 t_1 (1 + t_1 t_2))} {(-1 + t_1 t_2 t_3) (1 + t_1 t_2 t_3)}.
\end{align*}
The connection between two parametizations is given by
$$ t_1 = \frac{a_1}{r-1}, \quad
t_2 = \frac{-(1-r^2+a_1^2 t^2)}{2(t-1)a_1}, \quad
t_3 = \frac{2a_1 t(t-1)}{1-r^2+a_1^2 t^2}, $$
$$ a_1 =  \frac{2 t_1 (1 + t_1 t_2 (1 + t_2 t_3))} {(-1 + t_1 t_2 t_3) (1 + t_1 t_2 t_3)}, \quad
r = \frac{1+2t_1t_2+2t_1t_2^2t_3+t_2^2t_3^2t_1^2}{(-1+t_1t_2t_3)(1+t_1t_2t_3)},\quad
t = -t_2 t_3. $$

\section{New construction of families of Diophantine sextuples} \label{sec:sextuples}

Let $\{a_1,a_2,a_3,a_4\}$ and $\{a_1,a_2,a_3,a_5\}$ be regular Diophantine quadruples,
i.e. $a_4$ and $a_5$ are solutions of the quadratic equation
$$ (a_1+a_2-a_3-x)^2 - 4(a_1a_2+1)(a_3x+1) = 0. $$
We obtain that
\begin{align*}
a_4 &=\frac{-2(1-t_3+t_2t_3)(t_3t_1+1-t_1)(-t_2+1+t_1t_2)(-1+t_1t_2t_3)}{(1+t_1t_2t_3)^3}, \\
a_5 &=\frac{2(t_3+t_2t_3+1)(t_3t_1+t_1+1)(1+t_2+t_1t_2)(1+t_1t_2t_3)}{(-1+t_1t_2t_3)^3}.
\end{align*}
In order that $\{a_1,a_2,a_3,a_4,a_5\}$ be a rational Diophantine quintuple, it remains to satisfy the
condition that $a_4a_5+1$ is a perfect square. We obtain the condition that
\begin{align} \label{eq:pt1t2t3}
p(t_1,t_2,t_3)&=(-8t_2^3t_3^3-8t_3^2t_2^2-3t_3^4t_2^4+4t_2^2+4t_2^2t_3^4+4t_2^4t_3^2+8t_2^3t_3)t_1^4 \\
&\mbox{}+(8t_2^2t_3-16t_2t_3^2-8t_2^3t_3^2+8t_2-8t_2^3t_3^4-8t_2^4t_3^3-8t_2^2t_3^3+8t_3^4t_2)t_1^3 \nonumber \\
&\mbox{}+(-8t_3^2-8t_2^2-8t_2t_3-8t_2^3t_3^3-8t_2^2t_3^4+8t_3^3t_2+4t_3^4+4-18t_3^2t_2^2 \nonumber \\
&\mbox{}\,\,\,\,+4t_3^4t_2^4-8t_2^4t_3^2-16t_2^3t_3)t_1^2 \nonumber \\
&\mbox{}+(8t_2^4t_3^3-8t_2^2t_3-16t_2^2t_3^3-8t_2t_3^2-8t_2+8t_3^3+8t_2^3t_3^2-8t_3)t_1
\nonumber \\
&\mbox{}-3-8t_2t_3+4t_2^4t_3^2-8t_3^2t_2^2+4t_3^2+4t_2^2+8t_2^3t_3
\nonumber
\end{align}
is a perfect square. We compute the discriminant of $p(t_1,t_2,t_3)$ with the respect to $t_1$ and
factorize it. One of the factors is
$$ p_1(t_2,t_3)= 3+10t_2t_3-3t_3^2+3t_3^2t_2^2. $$
The condition $p_1(t_2,t_3)=0$ leads to $9t_3^2+16$ be a perfect square, say
$9t_3^2+16 = (3t_3+u)^2$. We get
\begin{align*}
t_3 &= \frac{16-u^2} {6u}, \\
t_2 &= \frac{u^2+10u+16} {(u-4)(u+4)}.
\end{align*}
Inserting this in (\ref{eq:pt1t2t3}), we obtain that $a_4a_5+1$ is a perfect square.
Thus, we obtained a two-parametric family (in parameters $t_1$ and $u$)
of rational Diophantine quintuples which contain two regular quadruples
(let us mention that a one-parametric family of rational Diophantine quintuples
with this property was constructed in \cite{duje-rocky1}).

Now we extend the nonregular quadruple $\{a_1,a_3,a_4,a_5\}$ to regular quintuples
$\{a_1,a_3,a_4,a_5,a_6\}$ and $\{a_1,a_3,a_4,a_5,a_7\}$, i.e. $a_6$ and $a_7$ are solutions
of the quadratic equation
$$ (a_1a_3a_4a_5x+2a_1a_3a_4+a_1+a_3+a_4-a_5-x)^2 -4(a_1a_3+1)(a_1a_4+1)(a_3a_4+1)(a_5x+1) = 0. $$
We will not use $a_7$ is our construction, so we give here only the value of $a_6$:
{\small
\begin{align*}
a_6 &=6(u+4)(u+8)(u+2)(u-4)(2t_1u^2+3u^2+20t_1u+12u+32t_1) \\
\mbox{}&\times (t_1u^2+10t_1u+16t_1-6u)(t_1u^2+10t_1u+16t_1+6u)(t_1u^2+10t_1u+16t_1-24-6u) \\
\mbox{}&\times (4096t_1^2+15360t_1^2u+15168t_1^2u^2+5920t_1^2u^3+948t_1^2u^4+60t_1^2u^5+t_1^2u^6-12288t_1u \\
\mbox{}&-7680t_1u^2+480t_1u^4+48t_1u^5-5184u^2-2592u^3-324u^4)^{-2}.
\end{align*}}%

The only missing condition in order that $\{a_1,a_2,a_3,a_4,a_5,a_6\}$ be a rational Diophantine
sextuple is that $a_2a_6+1$ is a perfect square. This condition leads to the quartic in
$t_1$ over $\mathbb{Q}(u)$:

{\footnotesize
\begin{align*}
& (u^{12}+120u^{11}+5496u^{10}+125600u^9+1639440u^8+13075200u^7+65656320u^6 \\
& \,\,\,\,+209203200u^5+419696640u^4+514457600u^3+360185856u^2+125829120u+16777216)t_1^4 \\
& +(24u^{12}++1296u^{11}+32256u^{10}+446208u^9+3461760u^8+13047552u^7-208760832u^5 \\
& \,\,\,\,-886210560u^4-1827667968u^3-2113929216u^2 -1358954496u-402653184)t_1^3 \\
& +(36u^{12}+1296u^{11}+18072u^{10}+48096u^9-1681632u^8-22516992u^7-127051776u^6 \\
& \,\,\,\,-360271872u^5-430497792u^4+197001216u^3+1184366592u^2 +1358954496u+603979776)t_1^2 \\
& +(-432u^{11}-15552u^{10}-259200u^9 -2267136u^8-9116928u^7+145870848u^5 \\
& \,\,\,\,+580386816u^4+1061683200u^3+1019215872u^2 +452984832u)t_1 \\
& +1296u^{10}+41472u^9+31643136u^6+670032u^8+6054912u^7+96878592u^5+171528192u^4 \\
& \,\,\,\, +169869312u^3+84934656u^2=z^2.
\end{align*}}%
Since this quartic has a $\mathbb{Q}(u)$-rational point at infinity,
it can be transformed by birational transformations
into an elliptic curve over $\mathbb{Q}(u)$
(the singular point at infinity on the quartic corresponds to the point at
infinity and an additional point $P_1$ on the elliptic curve).
The quartic have another $\mathbb{Q}(u)$-rational point
corresponding to $t_1=\frac{-3(u+4)u}{2(u^2+10u+16)}$.
It gives $a_6=0$, so it does not yield a rational Diophantine sextuples.
However, if we denote the corresponding point on the elliptic curve
by $P_2$, then the point $2P_2$ on the elliptic curve corresponds to the point with
$$ t_1= \frac{3(3u^4+40u^3+368u^2+1280u+1024)}
{4(u^2+10u+16)(u+20)u} $$
on the quartic, and by inserting this value, we obtain the parametric family of rational Diophantine sextuples


{\tiny
\begin{gather*}
 \left\{ 
\frac{-12u(u+4)(3u^4+8u^3+224u^2+576u+512)(3u^3+28u^2+256u+256)}
{(u+8)(u+2)(u-4)(3u^3+8u^2+144u+128)(3u^4+48u^3+528u^2+1280u+1024)}, \right. \\  
\frac{8u(u+20)(3u^5+8u^4+64u^3-640u^2-2304u-2048)(u+8)(u+2)}
{3(u+4)(u-4)(3u^3+8u^2+144u+128)(3u^4+48u^3+528u^2+1280u+1024)}, \\
\frac{2(u+4)(u-4)(39u^7+776u^6+8096u^5+48640u^4+226048u^3+587776u^2+770048u+393216)}
{3(u+8)(u+2)(3u^3+8u^2+144u+128)(3u^4+48u^3+528u^2+1280u+1024)}, \\
 \frac{-8(u^2+4u+32)(3u^3+14u^2-40u-64)(9u^3+8u^2+112u+384)(3u^4+48u^3+528u^2+1280u+1024)}
{3(u+8)(u+4)(u+2)(u-4)(3u^3+8u^2+144u+128)^3}, \\
\frac{4u(u+2)(17u^2+48u+48)(3u^5+8u^4-176u^3-2944u^2-9216u-8192)(3u^3+8u^2+144u+128)(u+8)^2}
{3(u+4)(u-4)(3u^4+48u^3+528u^2+1280u+1024)^3}, \\  
\left. \frac{12(u+2)(u-4)(5u+8)(u+4)(3u^2+8u+64)(3u^3+8u^2+144u+128)(3u^4+48u^3+528u^2+1280u+1024)}
{(u+8)(16384+69632u+64768u^2+22272u^3+3680u^4+576u^5+9u^6)^2}  \right\}
\end{gather*} }%
which satisfies the properties from Theorem \ref{tm1}.

\medskip

E.g. for $u=-1$ we get the rational Diophantine sextuple
$$ \left\{\frac{27900}{17479}, \,\, \frac{471352}{112365}, \,\, \frac{261770}{17479}, \,\, 
\frac{185535272}{419265}, \,\, \frac{63737828}{526368735}, \,\, \frac{79554420}{408480247} \right\}. $$
By taking other linear combinations of the points $P_1$ and $P_2$ we can also obtain
(more complicated) families of rational Diophantine sextuples.

\bigskip

{\bf Acknowledgements.}
The authors were supported by the Croatian Science Foundation under the project no.~IP-2018-01-1313.
The authors acknowledge support from the QuantiXLie Center of Excellence, a project
co-financed by the Croatian Government and European Union through the
European Regional Development Fund - the Competitiveness and Cohesion
Operational Programme (Grant KK.01.1.1.01.0004).
The authors acknowledge the usage of the supercomputing resources
of Division of Theoretical Physics at Ru\dj{}er Bo\v{s}kovi\'c Institute,
as well as the computing resources at Department of Mathematics, University of Zagreb
which were provided by Croatian Science Foundation grant HRZZ-9345


\begin{thebibliography}{99}

\bibitem{B-D}
A. Baker and H. Davenport, {\it The equations $3x^2 - 2 = y^2$ and $8x^2 - 7 = z^2$},
Quart. J. Math. Oxford Ser. (2) {\bf 20} (1969), 129--137.

\bibitem{BTF}
M. Bliznac Trebje\v{s}anin and A. Filipin,
{\it Nonexistence of $D(4)$-quintuples}, J. Number Theory {\bf 194} (2019), 170--217.

\bibitem{duje-quint}
A. Dujella,
{\it On Diophantine quintuples}, Acta Arith. {\bf 81} (1997), 69--79.

\bibitem{duje-rocky1}
A. Dujella,
{\it Diophantine triples and construction of high-rank elliptic curves over
$\mathbb{Q}$ with three non-trivial $2$-torsion points},
Rocky Mountain J. Math. {\bf 30} (2000), 157--164.

\bibitem{duje-irr}
A. Dujella, {\it Irregular Diophantine $m$-tuples and elliptic curves of high rank},
Proc. Japan Acad. Ser. A Math. Sci. {\bf 76} (2000), 66--67.

\bibitem{duje-crelle}
A. Dujella, {\it There are only finitely many Diophantine quintuples},
 J. Reine Angew. Math. {\bf 566} (2004), 183--214.

\bibitem{duje-sext}
A. Dujella, {\it Rational Diophantine sextuples with mixed signs},
Proc. Japan Acad. Ser. A Math. Sci. {\bf 85} (2009), 27--30.

\bibitem{Duje-notices}
A. Dujella, {\it What is ... a Diophantine $m$-tuple?}, Notices Amer. Math. Soc. {\bf 63} (2016), 772--774.

\bibitem{Duje-Matija}
A. Dujella and M. Kazalicki, {\it More on Diophantine sextuples}, in: Number Theory - Diophantine problems, uniform distribution and applications, Festschrift in honour of Robert F. Tichy's 60th birthday (C. Elsholtz, P. Grabner, Eds.), Springer-Verlag, Berlin, 2017, pp. 227--235.

\bibitem{DKMS}
A. Dujella, M. Kazalicki, M. Miki\'c and M. Szikszai, {\it There are infinitely many rational
Diophantine sextuples}, Int. Math. Res. Not. IMRN {\bf 2017 (2)} (2017), 490--508.

\bibitem{DKP-sext}
A. Dujella, M. Kazalicki and V. Petri\v{c}evi\'c,
{\it There are infinitely many rational Diophantine sextuples with square denominators}, preprint, 2019.

\bibitem{Gibbs1}
P. Gibbs,
{\it Some rational Diophantine sextuples}, Glas. Mat. Ser. III {\bf 41} (2006), 195--203.

\bibitem{G2}
P. Gibbs,
{\it Regular rational Diophantine sextuples}, preprint, 2016.

\bibitem{G3}
P. Gibbs,
{\it A survey of rational diophantine sextuples of low height}, preprint, 2016.

\bibitem{HTZ}
B. He, A. Togb\'e and V. Ziegler,
{\it There is no Diophantine quintuple}, Trans. Amer. Math. Soc., to appear.

\bibitem{Hea}
T. L. Heath, Diophantus of Alexandria. A Study in the History of Greek Algebra. Powell's Bookstore, Chicago; Martino Publishing, Mansfield Center, 2003.

\bibitem{luka}
L. Lasi\'c, personal communication, 2017.

\bibitem{P}
T. Piezas,
{\it Extending rational Diophantine triples to sextuples}, \\
{\tt http://mathoverflow.net/questions/233538/extending-rational-diophantine-triples-to-sextuples}

\bibitem{Stoll}
M. Stoll, {\it Diagonal genus 5 curves, elliptic curves over $\mathbb{Q}(t)$, and rational diophantine quintuples}, Acta Arith., to appear.



\end{thebibliography}
\end{document}